\documentclass[12pt]{amsart}
\usepackage{amssymb}
\usepackage{amsmath}

\newtheorem{thm}{Theorem}[section]
\newtheorem{lemma}[thm]{Lemma}

\newtheorem{theorem}[thm]{Theorem}
\newtheorem*{theorem*}{Theorem}
\newtheorem*{lemma*}{Lemma}

\newtheorem{proposition}[thm]{Proposition}

\newtheorem*{corollary*}{Corollary}

\theoremstyle{definition}

\numberwithin{equation}{section}
\newcommand{\st}[1]{\ensuremath{^{\scriptstyle \textrm{#1}}}}

\def    \Cinf   {C^\infty}

\newcommand{\CC}{{\mathbb C}}

\newcommand{\RR}{{\mathbb R}}
\newcommand{\ZZ}{{\mathbb Z}}
\newcommand{\bbC}{{\mathbb C}}
\newcommand{\bbR}{{\mathbb R}}

\newcommand{\fg}{{\mathfrak g}}

\newcommand{\Gr}{\mbox{Gr}}

\newcommand{\pr}{\mbox{pr}}

\newcommand{\trace}{\mbox{trace}}

\newcommand{\rmv}{\text{v}}
\newcommand{\rmw}{\text{w}}

\def \bar{\overline}

\numberwithin{equation}{section}

\theoremstyle{definition}
\theoremstyle{remark}

\newcommand{\F}{\mathcal{F}}

\newcommand{\C}{\mathbb{C}}

\newcommand{\T}{\mathbb{T}}

\def \Tp {T^\#}

\newcommand{\arabiclist}{
  \renewcommand{\theenumi}{\arabic{enumi}}%
  \renewcommand{\labelenumi}{\theenumi.}%
}

\newtheorem*{remarks*}{Remarks}

\begin{document}

\title{Geodesics on weighted projective spaces}
\author{V. Guillemin}\address{Department of Mathematics\\
Massachusetts Institute of Technology \\Cambridge, MA 02139 \\
USA}\thanks{V. Guillemin is supported in part by NSF grant
DMS-0408993.}\email{vwg@math.mit.edu}
\author{A. Uribe}\address{Department of Mathematics \\
University of Michigan \\ Ann Arbor \\ MI \\ 48109 \\ USA}
\thanks{A. Uribe is supported in part by NSF grant DMS-0401064.}
\email{uribe@umich.edu}
\author{Z. Wang}
\address{Department of Mathematics\\
Massachusetts Institute of Technology \\ Cambridge, MA  02139 \\
USA} \thanks{Z. Wang is supported in part by NSF grant
DMS-0408993.}\email{wangzq@math.mit.edu}

\date{}

\begin{abstract}
We study the inverse spectral problem for weighted projective spaces using
wave-trace methods.  We show that in many cases one can ``hear" 
the weights of a weighted projective space.  
\end{abstract}

\maketitle

\tableofcontents

\section{Introduction}
A weighted complex projective space is a quotient 
\[
X = S^{2d-1}/S^1,
\]
where $S^{2d-1}\subset\C^d$ is the unit sphere in $\C^d$ and the circle acts on the sphere as   
\[
e^{it}\cdot (z_1,\ldots ,z_d) = (e^{iN_1t}z_1,\ldots ,e^{iN_d t}z_d),
\]
where $N_1,\ldots N_d$ are positive integers, called the weights.
Although naturally a (singular) complex manifold, we will think of $X$ as a real Riemannian orbifold, of real dimension $2(d-1)$.

In \cite{ADFG}, Abreu, Dryden, Freitas and Godinho raise the interesting question
of whether the spectrum of the Laplace-Beltrami operator on $X$ determines the weights $N_1,\ldots , N_d$, and they show that this is indeed the case for $d=3$, if one knows the spectrum of the Laplacian
on zero and 1-forms.  Their methods center on the heat invariants and use equivariant
cohomology.

In a different direction, the wave-trace formula of \cite{DG} can be extended to compact Riemannian
orbifolds, see  \S 2 and \cite{SU}.  On the basis of this result one can expect that,
generically, the spectrum of the Laplacian on an orbifold determines the set of lengths (with multiplicities) of 
closed geodesics, together with the dimension of the manifold of closed 
geodesics of a given length. This is actually the case for weighted projective spaces. 
 Let's call the multiset of such pairs of data, (lengths of closed geodesics with multiplicities, and the dimension of their set), the weighted length spectrum.  From the wave-trace point of view, the inverse spectral problem of the previous 
paragraph is closely related to the problem:  {\em Does the weighted length
spectrum of $X$ determine the set of weights $\{N_j\}$?}

Geodesics on $X$ are projections of geodesics on $S^{2d-1}$ that are ``horizontal", meaning orthogonal to $S^1$ orbits.   Therefore all geodesics on $X$ are periodic with a common period, $2\pi$.
However, there can exist exceptional, shorter geodesics (they are ``exceptional" in the sense that the measure of their set in the space of all geodesics is zero).  They arise because of the non-trivial isotropy of some points on the sphere $S^{2d-1}$.

In this paper we show that, in many cases, the weighted length spectrum of 
$X$ corresponding to the exceptional geodesics
determines the set of sums of pairs of weights, $\{N_j+N_k, j<k\}$, with multiplicities.
This in turn raises the following problem:  {\em Is a set of $d$ positive numbers, $\{N_j\}$, determined by the multiset of the sums of the pairs $\{N_j+N_k, j<k \}$?} 
The answer is obviously ``no" for $d=2$ and obviously ``yes" for $d=3$.  It 
has been known for some time that the answer is ``yes" if $d$ is not a power of
two.  For convenience we present a proof of this result in an appendix, since
the only reference we found, \cite{MSU}, may not be readily available to the reader.
Conversely, Elizabeth Chen and Jeffrey Lagarias have shown recently
that, if $d$ is any power of two, then there exists sets of $d$ numbers that 
are not determined by the multiset of the sums of pairs.
(personal communication).

To briefly summarize our results, our calculations of the weighted length spectrum of $X$, together with the wave-trace formula for orbifolds, constitue a proof that for
$d$ not a power of two,
generically (with respect to the weights) one can ``hear" the weights of a weighted projective space.  

About the organization of the paper:  After reviewing the wave trace for Riemannian orbifolds in \S 2, we compute in \S 3 the lengths of the exceptional geodesics in the case $d=2$, which we feel makes clearer the general discussion in \S\S 4 and 5.  Our main inverse spectral result is Theorem \ref{thm:7.1}.

{\sc Acknowledgments:}  We are indebted to Jeffrey Lagarias for his help with the ``sums-of-pairs" inverse problem mentioned above.

\section{The wave trace for orbifolds}
\label{sec:1}
In this section we outline a proof of the ``wave trace formula"
for orbifolds, referring to \cite{SU} for details.
We recall that an orbifold can always be represented as the
quotient of a manifold by a group, i.e.,~as a quotient $X/G$
where $X$ is a manifold, $G$ a compact connected Lie group and
$\tau: G \times X \to X$ a locally free action of $G$ on $X$.
If in addition $X$ is equipped with a $G$-invariant Riemannian
metric the Laplace operator $\Delta_X : \Cinf (X) \to \Cinf (X)$
maps $\Cinf (X)^G$ into itself and via the identification, $\Cinf
(X)^G = \Cinf (X/G)$, defines a Laplace operator, $\Delta_{X/G}$ on $X/G$.
Assuming that $X$ is compact the wave trace formula for $X$
asserts that the expression
\begin{displaymath}
  e(t) = : \trace \exp it \sqrt{\Delta_X}
\end{displaymath}
is well-defined as a tempered distribution and that its singular
support is contained in (and modulo some slightly technical
clean intersection assumptions which we'll avoid going into at
this juncture) coincides with the period spectrum of
$X$.\footnote{A point, $T \in \RR$ is in the \emph{period
spectrum} of $X$ if there exists a closed geodesic, $\gamma$,
on $X$ of period $T$, i.e.,~$\gamma (t +T) = \gamma (t)$ for all
$t$.}  We will prove below that the analogue of this result is
also true for $X/G$ providing one defines carefully what one means
by a ``closed geodesic'' and by its ``period''.

Let $T^\#X$ be the punctured cotangent bundle of $X$.  From
$\tau$ one gets an action, $\tau^\#$, of $G$ on $\T^\# X$ and from
the symbol of $\Delta_X$ a Hamiltonian vector field, $v_H$, where
$H=\sigma (\Delta_X)^\frac12$.  Moreover $\tau^{\#} $ and the
flow generated by $v_H$ commute, so one has an action,
\begin{displaymath}
  (g,t) \to \tau^\#_g \exp t v_H
\end{displaymath}
of $G \times \RR$ on $T^\# X$.  The \emph{moment Lagrangian} for
this action (see \cite {We}) is the set of points
\begin{displaymath}
  (x,\xi,y,\eta ,t,\tau ,g, \gamma) \in T^\# (X \times X \times
  \RR \times  G)
\end{displaymath}
satisfying
\begin{equation}
  \label{eq:1.1}
  (y,\eta) = (\tau^\#_g \exp t v_H) (x,\xi)\, , \quad
  \tau = H (x,\xi)\, , \quad \gamma = \Phi (x,\xi)
\end{equation}
where $\Phi : \Tp X \to \fg^*$ is the moment map associated with
the action of $G$ on $\Tp X$.  The set (\ref{eq:1.1}) is also the
microsupport of the operator, $\tau^\#_g \exp it \sqrt{\Delta_X}$,
and we'll use this fact and functorial properties
of microsupports to compute the microsupport of the distribution
\begin{equation}
  \label{eq:1.2}
  \trace \exp it \sqrt{\Delta_{X/G}}
\end{equation}
Some observations:

\arabiclist
\begin{enumerate}

\item 
Modulo the identification, $\Cinf (X)^G = \Cinf (X/G)$ the expression
(\ref{eq:1.2}) can be written as an integral
  \begin{equation}
    \label{eq:1.3}
    \int_G \trace (\tau^\#_g \exp it \sqrt{\Delta_X})\, dg
  \end{equation}
where $dg$ is Haar measure.

\item 
Let $e(x,y,t,g)$ be the Schwartz kernel of the operator,
$\tau^\#_g \exp it \sqrt{\Delta_X}$.  Then (\ref{eq:1.3}) can
also be written as a double integral
\begin{equation}
  \label{eq:1.4}
  \int e (x,x,t,g)\, dg \, dx
\end{equation}
where $dx$ is the Riemannian volume form on $X$.

\item 
Let $\iota_\Delta : X \times \RR \times G \to X \times X \times
\RR \times G$ be the diagonal embedding and $\pi : X \times \RR
\times G \to \RR$ the projection, $(x,t,g) \to t$.  Then
(\ref{eq:1.4}) can be written more functorially as a
``pull-back'' followed by a ``push-forward'':
\begin{equation}
  \label{eq:1.5}
  \pi_* \iota^*_\Delta e (t) \, .
\end{equation}

\item 
Functorial properties of ``pull-backs'' tell us that if $(x,\zeta
,t,\tau ,g,\gamma)$ is in the microsupport of $\iota^*_\Delta e$
there exists a point $(x,\xi,y,\eta,t,\tau g,\gamma)$ in the
microsupport of $e$, i.e., in the set (\ref{eq:1.1}), such that
\begin{equation}
  \label{eq:1.6}
  x=y \hbox{   and   } \zeta = \xi -\eta
\end{equation}
and, by (\ref{eq:1.1}),
\begin{equation}
  \label{eq:1.7}
  (x,\eta) = (\tau^\#_g \exp t v_H) (x,\xi)\, .
\end{equation}

\item 
Functorial properties of ``push-forwards'' tell us that if
$(T,\tau)$ is in the microsupport of $\pi_* \iota^*_\Delta e$
there exists a point, $(x,\zeta,T,\tau ,g,\gamma)$ in the
microsupport of $\iota^*_\Delta e$, i.e., in the set
(\ref{eq:1.6})--(\ref{eq:1.7}) such that $\zeta = \gamma =0$.
Thus by (\ref{eq:1.6})--(\ref{eq:1.7}) we've proved:

\begin{theorem*}

If $T$ is in the singular support of the distribution,
(\ref{eq:1.5}) there exists an $(x,\xi,g) \in (\Tp X) \times G$
with
\begin{eqnarray}
  \label{eq:1.8}
(\Tp_g \exp T v_H) (x,\xi) = (x,\xi)\\
\noalign{\mbox{and}}\notag\\
\label{eq:1.9}
\Phi (x,\xi) =0 \, .
\end{eqnarray}
Moreover, modulo appropriate clean intersection assumptions, 
one can compute the leading term of the singularity of
(\ref{eq:1.5}) at $T$ in terms of the geometry of the
manifold of curves $(\Tp_g \exp t v_H) (x,\xi)$.

\end{theorem*}

\end{enumerate}

It remains to show that (\ref{eq:1.8}) and (\ref{eq:1.9}) can be
interpreted as saying that $T$ is in the period spectrum of
$X/G$.  Since the action of $G \times \RR$ on $\Tp X$ is locally
free the group of $(g,T)$'s satisfying (\ref{eq:1.8}) is a discrete
group.  Thus the projection, $G \times \RR \to \RR$ maps this
onto a discrete subgroup of $\RR$ and hence onto a lattice group
\begin{equation}
  \label{eq:1.10}
  \{ n T_0 \, , \, n \in \ZZ \} \, .
\end{equation}
Moreover since $G$ is compact the projection of this group onto
$G$ is a finite subgroup of order~$m$.  Thus by (\ref{eq:1.8})
\begin{displaymath}
  (\exp mT v_H) (x,\xi) = (x, \xi)
\end{displaymath}
showing that the curve
\begin{displaymath}
(\exp t v_H) (x,\xi) \, , \, -\infty < t< \infty
\end{displaymath}
is periodic of period $mT$ and hence its projection onto $X$ is a
closed geodesic, $\gamma (t)$, of period~$mT$.

Let's next unravel the meaning of (\ref{eq:1.9}).  This asserts
that the covector, $\xi \in \Tp_x$, is conormal to the orbit of
$G$ through $x$ and hence for the projected curve, $\gamma (t)$,
that $\frac{d\gamma}{dt} (0)$ is orthogonal at $x$ to the orbit
of $G$ through $x$.  Since $\Phi^{-1}(0)$ is $G \times \RR$
invariant, this holds as well at $\gamma (t)$, therefore the
implication of (\ref{eq:1.9}) is that at every point, $\gamma
(t)$, $-\infty <t<\infty$, the orbit of $G$ through $\gamma (t)$
is orthogonal to $\gamma (t)$.  In other words, $\gamma$ is
horizontal with respect to the quasi-fibration, $\pr : X \to X/G$,
and therefore its projection, $\pr \circ \gamma$, onto $X/G$ is a
geodesic on $X/G$.\footnote{Away from singular points
these projections are locally length-minimizing, and therefore worthy of the
name ``geodesics" on $X/G$.  Across singular strata they fail to be 
locally lenght-minimizing, but we will continue to call them geodesics, as no
locally length minimizing curves exist in that case.  Of
  course one has to check that this definition is independent of
  the presentation of $X/G$ as the quotient of $X$ by $G$, but
  this isn't hard to do.}  Finally note that if $T$ belongs to
the set (\ref{eq:1.10}) then by (\ref{eq:1.8}) there exists a $g
\in G$ such that
\begin{displaymath}
  (\tau^\#_g \exp T v_H) (x,\xi) = (x,\xi)
\end{displaymath}
and hence
\begin{displaymath}
  \gamma (0) = \tau_g \gamma (T) \, .
\end{displaymath}
Thus for the projection of $\gamma$ onto $X/G$
\begin{displaymath}
  \pr \circ \gamma (T) = \pr \circ \gamma (0) \, .
\end{displaymath}
In other words, $T$ is the period of the closed geodesic, $\pr \circ
\gamma (t)$, $0 \leq t \leq T$.  In addition, the calculation alluded to in the
previous theorem shows that in positive curvature there cannot be 
singularity cancellations from geodesics of a given length.
Therefore one has:

\begin{theorem*}

If $T \in \RR$ is in the singular support of the distribution
(\ref{eq:1.2}), there exists a closed geodesic on $X/G$ of period
  $T$.  Moreover, if appropriate clean intersection hypotheses
  are satisfied and $X$ has positive 
  sectional curvature, this result is an ``if and only if'' result.

\end{theorem*}

\section{Weighted projective lines}
\label{sec:2}

These spaces are obtained as quotients of the three sphere $S^3
\subseteq \CC^2$ by circle actions.  More explicitly, let $p$
and $q$ be positive integers with $1< p <q$ and g.c.d.$(p,q) =1$, and let
\begin{displaymath}
  \tau : S^1 \times S^3 \to S^3
\end{displaymath}
be the action $\tau (e^{i\theta}) (z_1,z_2) = (e^{ip}
z_1,e^{iq}z_2)$.  Then $X = S^3/S^1$ is an orbifold with
singularities at its north pole $[1,0]$ and south pole $[0,1]$.
Moreover we can define a Riemannian metric on the non-singular
part of $X$, as follows.  Let $v$ be the infinitesimal generator
of the action $\tau$.  Since this action is locally free $v
(p)\neq 0$ at all $p \in S^3$.  Let $\pi$ be the projection of
$S^3$ onto the quotient $X$.  Then, except at the north and
south poles, we can define the metric on $T_qX$, $q
= \pi (p)$, by identifying $T_q X$ with the horizontal part $H_p$
of $T_pS^3$, i.e.~with the space in $T_p$ orthogonal to $v_p$.  Since
the $S^1$ action is by isometrics this definition doesn't depend
on $p$.  Moreover if $\gamma : [0,2\pi] \to S^3$ is a geodesic
which is horizontal at $p = \gamma (0)$, i.e.~satisfies
\begin{displaymath}
 \Big\langle \frac{d\gamma}{dt}(0)\, , \, v_p \Big\rangle =0,
\end{displaymath}
then by parallel transport it's horizontal at $q=\gamma (t)$ so
it makes sense to talk about \emph{horizontal} geodesics in
$S^3$, and it's clear from our definition of the Riemann metric
on $X$ that their projections are geodesics on $X$.  In
particular $X$ is Zoll: all its geodesics are closed.

\subsection{Horizontal geodesics}
\label{sec:3}

Here's a more down-to-earth description of these horizontal
geodesics.  Every geodesic on $S^3$ is the intersection of $S^3$
with a two-dimensional subspace, $V$, of $\RR^4$.  Identify $\RR^4$
with $\RR^2 \oplus \RR^2$ and let $J: \RR^2 \to \RR^2$ be the map
$J (a,b) = (-b,a)$.  (In other words, via the identification
\begin{displaymath}
  \RR^2 \to \CC \qquad (x,y) \to x+ \sqrt{-1}\,y\, ,
\end{displaymath}
let $J$ be multiplication by $\sqrt{-1}$.)  Then the vector field
$v(\rmv)$ at $\rmv = (\rmv_1,\rmv_2) \in \RR^2 \oplus \RR^2$ is just
\begin{equation}
  \label{eq:3.1}
  \tilde{J} \rmv = (p J \rmv_1 , q J \rmv_2)\, .
\end{equation}
Thus if $V$ is a two-dimensional subspace of $\RR^4$ spanned by
$\rmv$ and $w$ with $|\rmv | = |w| =1$ and $\rmv \perp w$, and
$\gamma : [0,2\pi] \to S^3$ is the parametrized geodesic lying on
the circle $V \cap S^3$ with $\gamma (0) = \rmv$ and
$\displaystyle{\frac{d\gamma}{dt}(0) =w}$,
\begin{displaymath}
  \Big\langle v (\rmv),\frac{d\gamma}{dt}(0)\Big\rangle = 0
     \Leftrightarrow  \Big\langle \tilde{J}\rmv ,w \Big\rangle
     =0\, .
\end{displaymath}
We conclude:

\begin{theorem}
  \label{th:3.1}
Let $\tilde{\omega} \in \Lambda^2 (\RR^4)^*$ be the two-form
\begin{equation}
  \label{eq:3.2}
  \tilde{\omega} ({\rm v} ,w) = \langle \tilde{J}{\rm v} ,w \rangle
\end{equation}
and let $V \subseteq \RR^4$ be a two-dimensional subspace of
$\RR^4$.  Then the geodesic $V \cap S^3$ is horizontal iff $V$ is
Lagrangian with respect to the two-form $\tilde{\omega}$.

\end{theorem}

Note that if $\rmv = (\rmv_1,\rmv_2) \in \RR^2 \oplus \RR^2$ and
$w=(w_1,w_2) \in \RR^2 \oplus \RR^2$ then
\begin{equation}
  \label{eq:3.3}
  \tilde{\omega} (\rmv ,w) = p \omega (\rmv_1,w_1) + q \omega (\rmv_2,w_2),
\end{equation}
where
\begin{displaymath}
  \omega (\rmv_i,w_i) = \langle J \rmv_i , w_i \rangle \, .
\end{displaymath}
From this formula one easily deduces

\begin{theorem}
\label{th:3.2}
There are two kinds of Lagrangian subspaces of $\RR^4$:

\arabiclist
\begin{enumerate}
\item 
Subspaces of the form
\begin{equation}
  \label{eq:3.4}
  V_1 \oplus V_2 \subseteq \RR^2 \oplus \RR^2
\end{equation}
where each $V_i$ is a one-dimensional subspace of $\RR^2$.

\item 
Graphs of maps $A:\RR^2 \to \RR^2$ satisfying
\begin{equation}
  \label{eq:3.5}
  A^* \omega = - \frac{p}{q}\omega \, .
\end{equation}

\end{enumerate}

\end{theorem}

Note that if $V \subseteq \RR^4$ is a Lagrangian subspace of
type~1, $V \cap S^3$ is a horizontal geodesic whose projection
onto $X$ goes through the north and south poles and if $V$ is a
Lagrangian subspace of type~2, $V \cap S^3$ is a horizontal
geodesic whose projection onto $X$ doesn't go through
\emph{either} the north or south pole.  We'll call Lagrangians of
type~1 \emph{polar} Lagrangians and those of type~2
\emph{non-polar} Lagrangians.

\subsection{The action of $S^1$ on the Grassmannian of
  two-dimensional subspaces of $\RR^4$}
\label{sec:4}

To simplify slightly the exposition below, we'll assume that $p$
and $q$ are odd.  The linear action of $S^1$ on $\RR^2 \oplus
\RR^2$, defined by
\begin{equation}
  \label{eq:4.1}
  \tau (e^{it} ) = (\exp tpJ \, , \, \exp tq J),
\end{equation}
induces an action $\tau^\#$ of $S^1$ on the Grassmannian of
two-dimensional subspaces of $\RR^4$.  If $p$ and $q$ are both
odd this action isn't effective since $\tau (e^{i\pi}) =-I$.
However we do get an effective action of $S^1 / \{ \pm 1 \}$.
Next note the following:

\begin{lemma*}
  If $V$ is a two-dimensional subspace of $\RR^2 \oplus \RR^2$
  then either
  \begin{equation}
    \label{eq:4.2}
    V=V_1 \oplus V_2
  \end{equation}
where $V_1$ and $V_2$ are one-dimensional subspaces of $\RR^2$, or
\begin{equation}
  \label{eq:4.3}
  V={\rm graph}\,  A
\end{equation}
where $A: \RR^2 \to \RR^2$ is a map of the first summand of $\RR^2
\oplus \RR^2$ into the second, or
\begin{gather}
  \label{eq:4.3p}
  V= {\rm graph} \, A' \tag{\ref{eq:4.3}$^\prime$}
\end{gather}
%
where $A'$ is a map of the second summand of $\RR^2$ into the first.

\end{lemma*}

It's easy to see that $S^1 /\{ \pm 1 \}$ acts freely at points
(\ref{eq:4.2}) of $G_2 (\RR^4)$.  What about points of type
(\ref{eq:4.3}) or (\ref{eq:4.3p})?  These two cases are basically
the same so let's concentrate on (\ref{eq:4.3}).  Note that $\tau
(e^{it})$ maps the subspace (\ref{eq:4.3}) into itself iff
\begin{equation}
  \label{eq:4.4}
    (\exp tq J)A = A (\exp tp J),
\end{equation}
or, equivalently, iff
\begin{eqnarray}
  \label{eq:4.5}
  (\cos qt) A &=& (\cos pt) A\\
\noalign{\hbox{and}}\notag\\
\label{eq:4.6}
 (\sin qt) JA &=& (\sin pt) AJ\, .
\end{eqnarray}
Hence if we exclude the trivial case $A=0$ the element
$e^{it}$, $0 < t < \pi$, of $S^1$ acts trivially on the
Grassmannian at the point (\ref{eq:4.3}) iff
\begin{eqnarray}
  \label{eq:4.7}
  JA=AJ \hbox{   and   } e^{itp} = e^{itq}\\
\noalign{\hbox{or}}\notag\\
\label{eq:4.8}
  JA= -AJ \hbox{   and   }  e^{itp} = e^{-itq}\, .
\end{eqnarray}
But $e^{itp} = e^{itq}$ iff $e^{it(q-p)} =1$ and hence iff
\begin{equation}
  \label{eq:4.9}
  t=\frac{\pi k}{q-p} \hbox{   for some   } 0<k<q-p,
\end{equation}
and $ e^{itp} = e^{iqt}$ iff $e^{it (p+q)} =1$ and hence iff
\begin{equation}
  \label{eq:4.10}
  t=\frac{\pi k}{p+q} \hbox{   for some   } 0<k<q+p\, .
\end{equation}
Moreover if (\ref{eq:4.7}) or (\ref{eq:4.8}) holds and $A \neq 0$
then $A$ is invertible, and therefore graph~$A$ is also the graph of a map,
$A'$, from the second summand of $\RR^2 \oplus \RR^2$ to the
first.  Thus we've proved

\begin{theorem}
  \label{th:4.1}

Suppose $p$ and $q$ are odd.  Then from the action $\tau^\#$ one gets
an effective action of $S^1 / \{ \pm 1 \}$ which fixes the
subspaces $\RR^2 \oplus \{ 0 \}$ and $\{ 0 \} \oplus \RR^2$ and is
elsewhere locally free.  Moreover $\tau^\#$ is free except at
subspaces $V$ of type (\ref{eq:4.3}) with $AJ = \pm JA$.  If
$AJ =JA$ the isotropy group of $V$ in $S^1$ is $\ZZ_{q-p}$, and if
$AJ = -JA$ its isotropy group  is $\ZZ_{p+q}$.

\end{theorem}

If either $p$ or $q$ is even this result has to be slightly
modified.

\begin{theorem}
  \label{th:4.2}
Suppose that either $p$ or $q$ is even.  Then the action
$\tau^\#$ fixes the subspaces $\RR^2 \oplus \{ 0 \}$ and $\{ 0 \}
\oplus \RR^2$ and is elsewhere locally free except at subspaces
of type (\ref{eq:4.2}) and at subspaces of type (\ref{eq:4.3}),
with $AJ = \pm JA$.  If $V$ is of type (\ref{eq:4.2}) its
isotropy group is $\ZZ_2$, and if $V$ is of type (\ref{eq:4.3})
with $AJ = \pm JA$ its isotropy group is the same as in the
theorem above.

\end{theorem}

Since the action $\tau$ of $S^1$ on $\RR^4$ preserves the
symplectic form $\tilde{\omega}$, the induced action $\tau^\#$
on the Grassmannian of two-dimensional subspaces of $\RR^4$
preserves the set of Lagrangian subspaces.  Moreover, if $V = {\rm
graph}\,  A$ is Lagrangian then by (\ref{eq:3.5}) $\det A = -
\frac{p}{q}$.  Therefore, if $A$ is such that $AJ=JA$,
graph~$A$ is \emph{not} Lagrangian since  $\det A$ is positive.  
Hence from the theorems
above we get as a corollary:

\begin{theorem}
  \label{th:4.3}

If $p$ and $q$ are odd the action of $S^1 / \{ \pm 1 \}$ on the
Lagrangian Grassmannian is free except on the one-dimensional set
\begin{equation}
  \label{eq:4.11}
  \{ {\rm graph}\,  A_\lambda \, , \quad | \lambda |^2 =
      \frac{p}{q} \},
\end{equation}
where $A_\lambda : \CC \to \CC$ is the map $A_\lambda z =\lambda
\bar{z}$.  If either $p$ or $q$ is even the action of $S^1$ is
free except on the set of polar subspaces (\ref{eq:4.2}) (which
are stabilized by $\ZZ_2$) and on the set (\ref{eq:4.11}).
Moreover, no matter what the parity of $p$ and $q$, the
stabilizer group of points in the set (\ref{eq:4.11}) is $\ZZ_{p+q}$.

\end{theorem}

\subsection{Geodesics on the weighted projective line}
\label{sec:5}

By (\ref{eq:4.1}), $\tau^\# (e^{it})$ maps the subspace $V=
{\rm graph}\, A_\lambda$  to the subspace $V' = {\rm graph}\, A_{\lambda'}$, where
 $\lambda' = e^{it (p+q)}\lambda$.  Thus in particular $S^1$ acts
 transitively on the set of horizontal geodesics
 \begin{equation}
   \label{eq:5.1}
   \left\{ {\rm graph}\, A_\lambda \cap S^3 \, , \quad
     | \lambda |^2 = \frac{p}{q}\right\}.
 \end{equation}
Therefore these horizontal geodesics project to a single geodesic in
$X$.  We will call this the \emph{exceptional} geodesic.  Note
that this geodesic is non-polar, i.e.,~lies entirely in the
non-singular part of $X$.  Another key observation:  Suppose that
for some $0 < \theta < 2 \pi$, $\tau (e^{i\theta})$ maps a
Lagrangian subspace $V$ onto itself.  Then for every point
$\rmv$ on the geodesic $V \cap S^3$, $\rmv$ and its image,
$\tau (e^{i\theta})\rmv$, project onto the same point in $X$.  Now
let $w \in V \cap S^3$ be a vector in $V$ perpendicular to $\rmv$
and let
\begin{displaymath}
  \gamma (t) = \cos t \rmv + \sin t w \,  \quad
    0 \leq t \leq 2 \pi
\end{displaymath}
be the geodesic parametrization of the circle $V \cap S^3$, with
$\gamma (0) = \rmv$, $\frac{d}{dt} \gamma (0) =w$ and $\gamma
(t_0) = \tau (e^{i\theta} \rmv )$.  Since the action of $S^1$ on
the space of horizontal geodesics is locally free, $\theta =
\frac{2\pi}{k}$ for some $k \in \ZZ_+$ and $\gamma $ is the union
of $k$ geodesic segments
\begin{displaymath}
  \gamma_m (t) = \tau (e^{\frac{i 2\pi m}{k}}) \gamma (t) \, ,
     \quad 0 < t < t_0 \, .
\end{displaymath}
Thus $t_0 = \theta$ and each of these segments projects onto a
closed geodesic in $X$ of length $\frac{2 \pi}{k}$.  Hence from
Theorem~\ref{th:4.3} we conclude:

\begin{theorem}
  \label{th:5.1}

If $p$ and $q$ are odd the geodesics on the Zoll orbifold $X$
are all of length $\pi$ except for the exceptional geodesic
which is of length $\frac{2\pi }{ p+q}$.  If either $p$ or $q$ is even
the geodesics on $X$ are all of length $2\pi$ except for the
polar geodesics (which are of length $\pi$) and the exceptional
geodesic (which is again of length $\frac{2\pi }{ p+q}$).

\end{theorem}

Thus, coming back to the results of \S1 we've proved

\begin{theorem}
  \label{th:5.2}

From the spectrum of the Laplace operator on $X$ one can ``hear''
whether $p$ or $q$ is odd and can also hear their sum $p+q$.

\end{theorem}

\section{Circle actions on $\CC^d$}
\label{sec:6}

As in the previous sections we'll identify $\CC^d$ with $\RR^{2d}$, i.e.,~with the
sum
\begin{equation}
  \label{eq:6.1}
  \RR^2 \oplus \cdots \oplus \RR^2
\end{equation}
of $d$ copies of $\RR^2$, and we'll view the standard representation of $S^1$
on $\CC^d$ with weights $N_i$, $i=1,\ldots ,d$  as an
$\RR$-linear representation:
\begin{equation}
  \label{eq:6.2}
  \tau (e^{it}) = (\exp N_1 tJ , \ldots , \exp N_dtJ)\, .
\end{equation}
From (\ref{eq:6.2}) we get an action, $\tau^\#$, of $S^1$ on the
Grassmannian of two-dimensional subspaces of $\RR^{2d}$, and we'll
try to ascertain below what the stabilizer groups are for this
action.  Our approach to this problem will be an inductive one:
Assuming we know what these isotropy groups are for $\RR^{2d-2}$
we'll try to determine what they are for $\RR^{2d}$.  Henceforth
we'll assume that the $N_i$'s \emph{are greater than one and are
  pair-wise relatively prime}.

We first note that each of the summands of (\ref{eq:6.1}) is a
fixed point for the action $\tau^\#$, and that on the complement
of these $d$ fixed points this action is locally free.    If all the
$N_i$'s are odd then $\tau (e^{i\pi})=-I$ and therefore $\tau^\#
(e^{i\pi})$  is the identity map, and so, as in \S4, one gets from
$\tau^\#$ an effective action of $S^1 /\{ \pm 1\}$.  If, one the
other hand, $N_1$ is even, the action $\tau^\#$ itself is an
effective action; however, if $V$ is a two-dimensional subspace
of the sum
\begin{equation}
  \label{eq:6.3}
  \{ 0 \} \oplus \RR^2 \oplus \cdots \oplus \RR^2
\end{equation}
its stabilizer contains $\ZZ_2$.  Also, no matter what the parity
of the $N_i$'s
 is, if $V$ is a subspace which is spanned by $\rmv_1$ and
 $\rmv_2$ where $\rmv_1$ lies in (\ref{eq:6.3}) and $\rmv_2$ lies in
 the first summand of (\ref{eq:6.1}) then its stabilizer is
 \emph{exactly} $\ZZ_2$.

\noindent{\bf Definition.} We will call a two-dimensional
subspace, $V$, of $\RR^{2d}$ \emph{exceptional} if its stabilizer
is not the identity and is not $\ZZ_2$.

We'll now outline an inductive method for determining these
exceptional points in $\Gr_2 (\RR^{2d})$ and their stabilizers.
Let
\begin{equation}
  \label{eq:6.4}
  \RR^2 \oplus \{ 0 \} \oplus \cdots \oplus \{ 0 \}
\end{equation}
be the first summand of (\ref{eq:6.1}).

\begin{lemma}
  \label{lem:6.1}
If $V$ is exceptional then its projection onto (\ref{eq:6.4}) is
either zero or two-dimensional.

\end{lemma}

\begin{proof}
If the projection is one-dimensional 
then $V$ is spanned by vectors $\rmv_1$ and $\rmv_2$ where
$\rmv_1$ lies in (\ref{eq:6.4}) and $\rmv_2$ lies in
(\ref{eq:6.3}).  Therefore, as we noted above, the stabilizer of
$V$ is $\ZZ_2$ and $V$ is not exceptional.

\end{proof}


The same argument shows

\begin{lemma}
  \label{lem:6.2}

If $V$ is exceptional then its projection onto (\ref{eq:6.3}) is
either zero or two-dimensional.
\end{lemma}

From these results we conclude

\begin{lemma}
  \label{lem:6.3}

If $V$ is exceptional then either it has to be the space
(\ref{eq:6.4}), or be contained in the space (\ref{eq:6.3}), or be
the graph of a map $A$ of (\ref{eq:6.4}) into (\ref{eq:6.3}).
Moreover this map has to be injective.

\end{lemma}

If $V$ is the space (\ref{eq:6.4}) then it is a fixed point for
the action $\tau^\#$, and if $V$ is contained in (\ref{eq:6.3})
then it is an exceptional subspace of $\RR^{2d-2}$ and hence
determined by our induction procedure.  Let's focus therefore on
the case where $V = {\rm graph} \, A$, $A$ an injective map of
(\ref{eq:6.4}) into (\ref{eq:6.3}).  We can write this map as a
direct sum
\begin{displaymath}
  A= (A_2,\ldots ,A_d),
\end{displaymath}
where $A_k$ is a map of $\RR^2$ into the $k$\st{th} summand of
(\ref{eq:6.3}).  The following result is a consequence of (\ref{eq:4.7}) and
(\ref{eq:4.8}):

\begin{lemma}
  \label{lem:6.4}
Suppose $V$ is stabilized by $\tau^\# (e^{it})$
with $0<t<2\pi$ and $t \neq \pi$.  Then, for each $k\in\{\, 2,\ldots ,d\,\}$,
either

\begin{eqnarray}
 \rm{(i)}  & A_k &=0\,\, \mbox{or}\notag\\
\label{eq:6.5}
 \rm{(ii)} & JA_k &= A_kJ\,\, \mbox{and}\,\, e^{itN_1}=e^{itN_k}
      \,\, \mbox{or}\\
 \rm{(iii)} & JA_k &= -A_kJ \,\, \mbox{and}\,\,
      e^{itN_1}=e^{-itN_k}\, .\notag
\end{eqnarray}

\end{lemma}

\noindent{\emph{Some notation:}} Let $S_1 (A)$ be the set of $k$'s
for which (ii) holds and let $S_2(A)$ be the set of $k$'s for
which (iii) holds.  Also, if $k \in S_1 (A)$ let $\sigma_k =
|N_k-N_1|$ and if $k \in S_2 (A)$ let $\sigma_k (A) = N_1 + N_k$.

Putting together the lemmas above we get the main result of this
section.

\begin{theorem}
\label{th:6.5}
Let $  A= (A_2,\ldots ,A_d)$ be a above, and
let $S(A) = S_1 (A) \cup S_2 (A)$.  Then the graph of $A$ is exceptional
if the components of $A$ satisfy (\ref{eq:6.5}) and the greatest common
divisor of the set of numbers
\begin{equation}
  \label{eq:6.6}
  \{ \sigma_k (A) \, , \quad k \in S (A) \}
\end{equation}
is greater than $2$.

\end{theorem}

Thus one can determine the isotropy groups of the exceptional
points of $\Gr_2 (\RR^{2d})$ by the following procedure:

\arabiclist
\begin{enumerate}
\item 
Fix two disjoint subsets $S_1$ and $S_2$ of the set $\{ 2 ,\ldots
,d \}$.

\item 
For each $k \in S_1$ let $\sigma_k = |N_k - N_1 |$.

\item 
For each $k \in S_2$ let $\sigma_k = N_k + N_1 $.

\item 
Let $\sigma (S_1,S_2)$ be the greatest common divisor of the set
of numbers
\begin{displaymath}
\{\,  \sigma_k \, , \quad k \in S_1 \cup S_2\,\}\, .
\end{displaymath}

\end{enumerate}

Then if $N=\sigma (S_1,S_2)$ is greater than $2$, $\ZZ_N$ is the
isotropy group of an exceptional point, $V= {\rm graph}\, A$, of
$\Gr_2 (\RR^{2d})$.

This observation gives us, by induction, the result below.

\begin{theorem}
  \label{th:6.6}

  For $1 \leq r \leq d-1$ let $S_1$ and $S_2$ be disjoint subsets
  of the set $\{ r+1 , \ldots, d \}$, and let $\sigma_k =
  |N_k-N_r|$ if $k$ is in $S_1$ and $\sigma_k= N_k + N_r$ if $k$ is in
  $S_2$.  Then, if the greatest common divisor $N=\sigma
  (S_1,S_2)$ of the set of numbers $\{\,\sigma_k$, $k \in S_1 \cup
  S_2\,\}$ is greater than $2$, $\ZZ_N$ is the isotropy group of an
  exceptional point of $\Gr_2 (\RR^d)$.  

\end{theorem}

\section{Geodesics on weighted $(d-1)$-dimensional projective
spaces}
\label{sec:7}

As mentioned in the introduction, this space is the quotient of $S^{2d-1}$ by the
circle action (\ref{eq:6.2}), i.e.,~it is the space
\begin{equation}
  \label{eq:7.1}
  X= S^{2d-1} /S^1 \, .
\end{equation}
(The title of this section is a bit of a misnomer since $X$, as a
real $\Cinf$ Riemannian orbifold, is ($2d-2$)-dimensional.) There are ${d}\choose {2}$ totally-geodesic embedded weighted
projective lines $X_{i,j} \ 1\leq i<j\leq d$
in $X$, corresponding to the direct sum of two complex
coordinate lines in $\bbC^d$.  As we showed in \S 3, each weighted
projective line has a unique exceptional geodesic $\gamma_{i,j}$,
and it has length 
\[
\ell_{i,j} = \frac{2\pi}{N_i+N_j}
\]
(assuming, as we always do in this paper, that the weights are pairwise 
relatively prime).  In this section we 
investigate conditions on the weights which guarantee that the lengths
$\ell_{i,j},\ 1\leq i<j\leq d$, are the {\em shortest} in the length spectrum of $X$.
As a corollary of the discussion of \S 2 in such cases one can ``hear" the 
multiset of the $\ell_{i,j}$, and therefore the multiset of pair-wise sums $\{ N_i+N_j,
1\leq i<j\leq d\}$.  This analysis requires that we have a description of the lengths of
the exceptional geodesics on $X$ (other than the $\gamma_{i,j}$).

Let
$\tilde{\omega} \in \Lambda^2 (\RR^{2d})^*$ be the symplectic form
\begin{displaymath}
  \tilde{\omega} (\rmv ,\rmw) = \sum^d_{i=1} N_i \langle
     J\rmv_i,\rmw_i \rangle
\end{displaymath}
where $\rmv_i$ and $\rmw_i$ are the projections of $\rmv$ and $\rmw$
onto the $i$\st{th} summand of (\ref{eq:6.1}).  The geodesics
on the sphere $S^{2d-1}$ are just the intersections $V \cap
S^{2d-1}$, where $V$ is a two-dimensional subspace of $\RR^{2d}$,
and one can show by exactly the same argument as in \S3 that $V
\cap S^{2d-1}$ is \emph{horizontal} if and only if $\iota^*_V
\tilde{\omega}=0$, i.e.,~if and only if $V$ is an isotropic
two-dimensional subspace of $\RR^{2d}$.  Given such a geodesic, the
projection $S^{2d-1} \to X$ maps $V \cap S^{2d-1}$ onto a closed
geodesic of $X$ of period $\frac{2\pi}{m}$, where $\ZZ_m$ is the
stabilizer of $V$ with respect to the action $\tau^\#$ of $S^1$
on $\Gr_2 (\RR^{2d})$.
In particular, suppose $V= {\rm graph} (A)$ as in the previous section.  
Then $V$ is isotropic if and only if
\begin{equation}
  \label{eq:7.2}
  N_1 + \sum^d_{m=2} N_m \det A_m =0,
\end{equation}
and if this is the case there exists $m$ such that $\det A_m$ is negative.  
Therefore, if
$V$ is exceptional and isotropic, (\ref{eq:6.5})(iii) must hold for some $m$,
or, equivalently, $S_2(A)\not=\emptyset$.
It follows that the stabilizer
group of $V$ has to be contained in $\ZZ_{N_1 + N_m}$ for
some~$m$. 

We summarize:
\begin{proposition}\label{k}
The lengths of the exceptional geodesics on $X$ are of the form $\ell_k = {2\pi}/k$,
where $k$ is the greatest common divisor of a set of numbers of the form
\[
\F_{S_1,\, S_2} = \{\, \sigma_m\;;\; m\in S_1\cup S_2 \},
\] 
and where $S_1,\, S_2$ are sets such that
\[
S_1\cup S_2\subset \{r+1,\ldots d\},\ S_1\cap S_2= \emptyset,\ 
S_2\not=\emptyset\,
\]
and such that g.c.d.$(\F_{S_1,\, S_2}) > 2$.  Here
\[
\sigma_m = 
\bigl\{
\begin{array}{c}
|N_m-N_r|\quad  \text{if}\ m\in S_1\\
N_m+N_r\quad  \text{if}\ m\in S_2.
\end{array}
\bigr.
\]
\end{proposition}

Note that if one takes $S_1=\emptyset$ and $S_2 = \{m\}$ (a singleton), then 
$k=N_r+N_m$ gives rise to the length of $\gamma_{r,m}$.

\bigskip
In low dimensions we have the following results:
\begin{theorem}
  \label{thm:7.1}
  Suppose all the $N_i$'s are pair-wise relatively prime.  Then:\\
  $\mathrm{(1)}$ If $d=3$, one can ``hear" the weights
                 $N_1, N_2$ and $N_3$. \\
  $\mathrm{(2)}$ If $d=4$, one can ``hear" all the pairwise
                 sums $\{\,N_i+N_j,\; i<j\,\}$, which will determine at most
                 two different choices of $N_i$'s. 
\end{theorem}

\subsection{The case $d=3$}
\label{sec:7.1}

The exceptional
geodesics other than the $\gamma_{i,j}$, if they exist, have to be associated with
isotropics of the form $V = {\rm graph} \, A$, where $A =
(A_2,A_3)$ and each of the $A_i$'s is bijective.  If we order
the $N_i$'s so that $N_1 < N_2 < N_3$, these geodesics will have lengths
of the form $2\pi /k$ where the possibilities for $k$ are
\begin{eqnarray}
  \label{eq:7.3}
  k &=& \text{g.c.d} (N_1 + N_3 , N_2 - N_1)\\
\noalign{\hbox{or}}\notag \\
  \label{eq:7.4}
  k &=& \text{g.c.d} (N_1 + N_2 , N_3 - N_1)\\
\noalign{\hbox{or}}\notag \\
  \label{eq:7.5}
  k &=& \text{g.c.d} (N_1 + N_2 , N_1 + N_3)\, .
\end{eqnarray}
(These possibilities correspond to the set $S_2$ of Theorem \ref{th:6.5}
being equal to $\{3\}$, $\{2\}$ and $\{2, 3\}$, respectively.)
If $k$ is given by (\ref{eq:7.3}) then $k < N_2<N_1+N_2$, and the
geodesic corresponding to $k$ is strictly longer  than the
geodesics $\gamma_1$, $\gamma_2$ and $\gamma_3$.  If $k$ is given
by (\ref{eq:7.4}) then $k<N_1+N_3$ and therefore the geodesic
is strictly longer than $\gamma_2$ and $\gamma_3$.
The geodesic will also be longer than $\gamma_1$ unless $N_1 +
N_2$ divides $N_3 -N_1$, or, equivalently (since $N_2 + N_3 = N_3 -
N_1 + N_1 + N_2$), unless $N_1 + N_2 $ divides $N_2 + N_3$.
Similarly, if $k$ is given by (\ref{eq:7.5}) the geodesic will be
strictly longer unless
$N_1 + N_2$ divides $N_1 + N_3$.  Thus if neither of these two
worst case scenarios occurs, $\gamma_1$, $\gamma_2$ and $\gamma_3$
are the three shortest geodesics on $X$. Moreover, they lie on the
non-singular part of $X$ and are isolated and non-degenerate.
Hence by \cite{DG}, $\gamma_i$ makes a contribution of the form
\begin{equation}
  \label{eq:7.6}
  c_i (t-T_i)^{-1}_+ + \cdots
\end{equation}
to the wave trace with $c_i \neq 0$ where $T_i$ is the period of
$\gamma_i$ and hence one can ``hear'' the $T_i$'s, i.e.,~one can
hear $N_1 + N_2$, $N_1 + N_3$ and $N_2 + N_3$ and hence one can
hear $N_1,N_2$ and $N_3$.

Let's next examine the first worst case scenario:  $N_1 + N_2$ divides
$N_1 + N_3$.  Then $V= {\rm graph}\, A$ is isotropic and stabilized
by $\ZZ_{N_1 + N_2}$ if and only if $A= (A_{\lambda_2},
A_{\lambda_3}) $ with
\begin{equation}
  \label{eq:7.7}
  N_1 - N_2 |\lambda_2 |^2 - N_3 |\lambda_3 |^2 = 0
\end{equation}
(see (\ref{eq:4.11})).  Thus the set of geodesics stabilized by
$\ZZ_{N_1 + N_2} $ consists of the isolated geodesic $\gamma_1$
on the projective line $X_1$, and completely disjoint from it the
2-parameter family of geodesics lying on the quotient by $S^1$ of
the three sphere (\ref{eq:7.7}).  In this case the wave trace
contains contributions at $T_1 = 2\pi / (N_1 + N_2)$ from both
$\gamma_1$ and from the geodesics (\ref{eq:7.7}).  However, since
the geodesics (\ref{eq:7.7}) also lie on the non-singular part of
$X$ the clean intersection techniques of \cite{DG} show that their
contribution to the wave trace is a singularity of the form
\begin{equation}
  \label{eq:7.8}
  c (t-T_1)^{-3} + \cdots
\end{equation}
which dominates the singularity (\ref{eq:7.6}), so again one can
hear $N_1 + N_2$.

A similar situation occurs if $N_1 + N_2$ divides $N_2 + N_3$. In
this case one gets again a 2-parameter family of geodesics on $X$
corresponding to $A$'s of the form $A= (A_{\lambda_2},
A_{\lambda_3})$ where
\begin{equation}
  \label{eq:7.9}
  N_2 | \lambda_2 |^2 = N_1 + N_3 |\lambda_3 |^2\, .
\end{equation}
The projection of this set onto the space of geodesics on $X$
isn't closed, but its closure  consists of the geodesics in this
set plus the geodesic $\gamma_3$ (which is stabilized by $\ZZ_{N_2
+ N_3}$ and hence by $\ZZ_{N_1 + N_2}$).  Thus the clean
intersection techniques of \cite{DG} show that the contribution of
this family to the wave trace is of the form (\ref{eq:7.8}), which
again dominates the contribution (\ref{eq:7.6}) coming from
$\gamma_1$.  Thus, in both these two worst cases scenarios one can
hear $N_1 + N_2$  and therefore all three $N_i$.

\subsection{The case $d=4$}
\label{sec:7.2}

Without loss of generality, we suppose $N_1 < N_2 < N_3 <N_4$. 
There are now ${4\choose2}=6$ embedded totally 
geodesic weighted projective lines, each one with a unique exceptional geodesic 
$\gamma_{i,j}$,
of length $\frac{2\pi}{N_i+N_j}$ (these are the ``desirable" lengths that we wish
to ``hear").   Note that $\gamma_{4,3}$ is the shortest geodesic, $\gamma_{4,2}$
the second shortest, $\gamma_{2,1}$ the longest and $\gamma_{3,1}$ the 
second longest of these geodesics, and in general we can't order the lengths of
$\gamma_{4,1}$ with $\gamma_{3,2}$ (in fact they could be equal, e.g.:
$17+3 = 13+7$).

The lengths of undesirable
exceptional geodesics are of the form
$L_k = \frac{2\pi}{k}$, $k$ a positive integer.  
We want to bound the possible values of $L_k$ from below, to determine which
desirable lengths are spectrally determined.

According to our inductive procedure for finding the possible $k$'s, we should first consider possible two-dimensional planes contained in 
\[
\{ 0\}\oplus \bbR^2\oplus \bbR^2\oplus \bbR^2.
\]
These correspond to the exceptional 
geodesics of a $d=3$ weighted projective space, with 
weights $N_2,\, N_3,\, N_4$.  As we have seen, the possibilities 
for $k$ are the greatest common divisor of the following sets:
\begin{equation}\label{uno}
\{N_4+N_2,\,N_3-N_2\},\quad 
\{N_3+N_2,\,N_4-N_2\},\quad 
\{N_3+N_2,\,N_4+N_2\}.
\end{equation}
The remaining exceptional geodesics are graphs of maps 
$A = (A_2,\,A_3,\,A_4)$.  The possible values of $k$ are the greatest common
divisors of the following seven sets:
\begin{equation}\label{dos}
\begin{array}{cc}
\{N_2+N_1,\,N_3+N_1\}, &\{N_2+N_1,\,N_3+N_1,\, N_4-N_1\}, \\
\{N_2+N_1,\,N_4+N_1\}, &\{N_2+N_1,\,N_4+N_1,\, N_3-N_1\}, \\
\{N_3+N_1,\,N_4+N_1\}, &\{N_3+N_1,\,N_4+N_1,\, N_2-N_1\}
\end{array}
\end{equation}
and
\begin{equation}\label{tres}
\{N_2+N_1,\,N_3+N_1,\, N_4+N_1\}.
\end{equation}

Using only that the g.c.d. of any set of numbers is 
less than or equal to any element in the set, one can check that in all
cases $k\leq N_3+N_2$, so we can always hear the lengths of the 
two shortest $\gamma_{i,j}$ and determine $N_4+N_3$ and $N_4+N_2$.
In addition, we claim that we can always hear $N_4+N_1$, because the
g.c.d.'s of all the sets listed above are all less than $N_4+N_1$.
The only set for which this is not immediately obvious is $\{N_3 + N_2,\, N_4+N_2\}$;
however, note that
\[
\text{g.c.d}\{N_3 + N_2,\, N_4+N_2\} = 
\text{g.c.d}\{N_3 + N_2,\, N_4-N_3\} < N_4 < N_4+N_1.
\]

One can also check that the g.c.d of the sets (\ref{dos})  satisfies
$k\leq N_3+N_1$, and
the g.c.d of the sets (\ref{tres}) $k\leq N_1+N_2$.
Therefore the lengths of the corresponding geodesics are greater than
the length of $\gamma_{3,2}$. If the g.c.d.'s of the sets (\ref{uno}) is less
than $N_3+N_2$, then all undesirable exceptional geodesics are longer
than $\gamma_{3,2}$ as well.   However, by our analysis of the ``worst case scenario"
when one of the g.c.d.'s of the sets (\ref{uno}) equals $N_3+N_2$, we can
conclude that we can always hear  $N_3+N_2$.

In conclusion, the multiset $\{ N_2+N_3, N_1+N_4, N_2 +N_4,
N_3+N_4\}$ is always spectrally determined.
In this set we can identify the individual sums $ N_3 +N_4$ and
$N_2+N_4$ as the largest and second largest elements,
respectively.  Moreover, although we are not able to distinguish between the
two remaining elements in the set which one is
$N_2+N_3$ and which one is $N_1+N_4$, we can form their sum,
\[
N_1+N_2+N_3+N_4, 
\]
and from this determine $N_1+N_2$ and $N_1+N_3$,
(by subtracting the known quantities $ N_3 +N_4$ and
$N_2+N_4$).  This shows that we can always hear the set
of all six sums of pairs $\{N_i+N_j, i<j\}$.

To reconstruct the possible weights giving rise to a given multiset of 
sums of pairs, proceed as follows:  $N_1+N_2$ and $N_1+N_3$
are the smallest and second smallest elements of the set,
respectively.  However, $N_2+N_3$ may be the third or fourth
element in the set if the set has six elements (the set may have only
five elements, exactly when $N_2+N_3 = N_4+N_1$, in which case this sum 
is the middle element and the reconstruction of the weights is unique).
Each one of those possibilities for $N_2+N_3$ gives a unique
solution for $N_1,\, N_2,\, N_3$, and the remaining elements
in the set uniquely determine the fourth weight.

It is clear that picking $N_2+N_3$ to be the third or the fourth element
from a set of six distinct pairwise sums gives different solutions $\{N_j\}$.
However, the two solutions ``usually" do not both consist of four pair-wise 
relatively prime numbers, although this can happen.  For example,
\begin{equation*}
   \{25, 29, 41, 61\} \mbox{\ and\ } \{17, 37, 49, 53\}
\end{equation*}
have the same set of pairwise sums.

\subsection{Remarks for $d$ general}

In general, if the ${d}\choose{2}$ shortest geodesics of $X$ are the
$\gamma_{i,j}$ then one can ``hear" the weights if $d$ is not a power
of two.  A way to ensure that the $\gamma_{i,j}$ are the 
${d}\choose{2}$ shortest geodesics is the following:

\begin{theorem}
 If $N_1<N_2<\ldots <N_d$, 
 $N_d \le 2 N_1$ and the list $N_i+N_j\,,\,1\leq i< j\leq d,$ contains no
repeated elements, then the  $\gamma_{i,j}$ are the ${d}\choose{2}$ shortest geodesics,
(and therefore, by the result in the appendix, the set of weights is spectrally
determined if $d$ is not a power of two).
\end{theorem}
\begin{proof}
Let $\frac{2\pi}{k}$ be an undesirable length of an exceptional geodesic.
How short can the exceptional geodesic be or, equivalently, how large can $k$ be?
Using Proposition \ref{k}, it is not hard to see that 
\[
k\leq \text{g.c.d}\{N_i+N_j,\, N_l+N_m\}
\]
for some $i<j$, $l<m$ and $(i,j)\not= (l,m)$.  Thus it suffices to show that such g.c.d.'s are
less than $N_1+N_2$.

First, the condition  $N_d \le 2 N_1$ implies that 
$N_l+N_m < 2(N_i+N_j)$.  Therefore 
\[
\text{g.c.d}\{N_i+N_j,\, N_l+N_m\}\not= N_i+N_j
\]
(since $a=\text{g.c.d}\{a,b\}\,\wedge\, b<2a\Rightarrow b=a$, and we assumed 
that the list $N_i+N_j\,,\,1\leq i< j\leq d,$ contains no
repeated elements).  It follows that $k\leq (N_i+N_j)/2 < N_1+N_2$.
\end{proof}

\begin{appendix}
\section{}
We present in this appendix a proof of the following:
\begin{theorem}
If $d$ is not a power of two, a set of $d$ real numbers
$\{N_1,\ldots ,N_d\}$ is determined by the multiset of sums of pairs
$\{N_j+N_k\ :\ j < k\}$.
\end{theorem}
\begin{proof}
The proof is based on the identity
\begin{equation*}
\sum_{1\le i < j \le d}(N_i+N_j)^k = (d-2^k)\sum_{i=1}^d N_i^k
+ \frac 12 \sum_{m=1}^{k-1}{k \choose m}(\sum_{i=1}^d N_i^{k-m})(\sum_{i=1}^d N_i^m).
\end{equation*}

This shows that if $d$ is not a power of two then one can find $\sum N_i^k$
from $\sum_{1\leq i<j\leq d}(N_i+N_j)^k$ and the sums $\sum N_i^m$ with $m<k$.
Proceeding inductively (starting from $\sum_{1\leq i<j\leq d}(N_i+N_j) = (d-1)\sum N_i$),
we can conclude that the elementary symmetric functions on the sums of pairs
determine the elementary symmetric functions on the $N_j$.
\end{proof}

Finally, we remark that one can generalize the previous discussion of the $d-2$ case
to obtain an algorithm for finding all possible sets of $N_i$'s, given the multiset of
their pair-wise sums.  This algorithm shows that, if 
$d$ is a power of $2$ and $d>2$, then the multiset $\{N_i+N_j\}$
will determine at most $d-2$ different choices of $N_i$'s.  However, 
it seems very hard to determine what the sharp bound of possible choices of $N_i$'s is, for a given $d$.

\end{appendix}

\end{document}